\documentclass{article}

\usepackage{arxiv}
\usepackage[utf8]{inputenc} 
\usepackage[T1]{fontenc}   
\usepackage{booktabs}       
\usepackage{amsfonts}
\usepackage{amsmath}
\usepackage{nicefrac}       
\usepackage{microtype}    
\usepackage{graphicx}
\usepackage{doi}
\usepackage{array,collcell} 
\usepackage{makecell}
\usepackage{mathtools}
\usepackage{nomencl}
\usepackage{tikz}
\usetikzlibrary{positioning}
\usepackage{pgfplots}
\pgfplotsset{compat = newest}
\usepackage{amssymb}
\usepackage{algorithmic}
\usepackage{textcomp}
\usepackage{cleveref}    
\usepackage{ntheorem}

\usepgfplotslibrary{groupplots}
\definecolor{butter1}{rgb}{0.988,0.914,0.310}
\definecolor{chocolate1}{rgb}{0.914,0.725,0.431}
\definecolor{chameleon1}{rgb}{0.541,0.886,0.204}
\definecolor{skyblue1}{rgb}{0.447,0.624,0.812}
\definecolor{plum1}{rgb}{0.678,0.498,0.659}
\definecolor{scarletred1}{rgb}{0.937,0.161,0.161}

\newtheorem{theorem}{Theorem}
\newtheorem{prop}{Proposition}

\title{Speeding up  backpropagation  of gradients through the Kalman filter via closed-form expressions}

\renewcommand{\*}[1]{\mathnormal{#1}}
\renewcommand{\!}[1]{\mathcal{#1}}
\renewcommand{\$}{\partial}

\newcommand{\fracd}[2]{\frac{\$#1}{\$#2}}

\newcommand{\fracdlb}[1]{\fracd{\!{L}}{#1}}
\newcommand{\dv}[2]{{\frac{\partial #1}{\partial #2}}}

\newcommand{\ipriorn}{_{n|n-1}}
\newcommand{\ipostn}{_{n|n}}
\newcommand{\ipostnn}{_{n-1|n-1}}

\newcommand{\xhpn}{\*{\hat{x}}\ipostn}
\newcommand{\xhpnn}{\*{\hat{x}}\ipostnn}
\newcommand{\xhmn}{\*{\hat{x}}\ipriorn}

\newcommand{\Ppn}{\*P\ipostn}
\newcommand{\Ppnn}{\*P\ipostnn}
\newcommand{\Pmn}{\*P\ipriorn}

\newcommand{\Id}{\*I}
\newcommand{\IKH}{\*I - \Kn \Hn}
\newcommand{\Kn}{\*K_n}
\newcommand{\Fn}{\*F_n}

\newcommand{\Sn}{\*S_n}
\newcommand{\Sninv}{\*S_n^{-1}}
\newcommand{\Hn}{\*H_n}
\newcommand{\Rn}{\*R_n}
\newcommand{\Qn}{\*Q_n}
\newcommand{\yn}{\*y_n}
\newcommand{\zn}{\*z_n}

\newcommand{\Lc}{\mathcal{L}}
\newcommand{\lnp}{l\ipostn}
\newcommand{\lnm}{l\ipriorn}
\newcommand{\Linc}{\Lc_{NLL}}

\newcommand{\squad}{\mkern8mu} 
 
\newcommand\AddLabel[1]{
  \refstepcounter{equation}
  (\theequation)
  \label{#1}
}
\newcolumntype{m}{>{\hfil$\displaystyle}X<{$\hfil}} 
\newcolumntype{n}{>{\collectcell\AddLabel}r<{\endcollectcell}}

\graphicspath{{images/}}

\author{ 
        {Colin Parellier} \\
	MINES Paris, PSL Research University\\ 
        Centre for Robotics,\\ 
        60 Bd Saint-Michel, 75006 Paris, France. \\ 
        Safran Tech, Groupe Safran, \\ 
        Rue des Jeunes Bois - Chateaufort, \\
        78772, Magny Les Hameaux, France.\\
	\texttt{colin.parellier@minesparis.psl.eu} \\
         \And
        {Silvere Bonnabel} \\
	MINES Paris, PSL Research University\\ 
        Centre for Robotics,\\ 
        60 Bd Saint-Michel, 75006 Paris, France. \\ 
	\texttt{silvere.bonnabel@minesparis.psl.eu} \\
	\And
        {Axel Barrau} \\
	OFFROAD\\ 
        5ter rue Parmentier, 94140 Alfortville.\\
}

\renewcommand{\shorttitle}{Speeding up  backpropagation of gradients through the KF via closed-form expressions}

\hypersetup{
pdftitle={Speeding up  backpropagation  of gradients through the Kalman filter via closed-form expressions},
pdfauthor={C. Parellier, A. Barrau, and  S. Bonnabel},
pdfkeywords={Kalman filter, maximum likelihood, backpropagation, matrix derivative, sensitivity tuning},
}

\begin{document}
\maketitle

\begin{abstract}
In this paper we provide  novel closed-form  expressions enabling differentiation of any  scalar function of the Kalman filter's outputs with respect to all its tuning parameters and to the measurements.  The approach differs from the previous well-known sensitivity equations in that it is based on a backward (matrix) gradient calculation,  that leads to drastic reductions  of the  overall computational cost. It is our hope that  practitioners  seeking numerical efficiency  and reliability will benefit from the concise  and exact equations  derived in this paper and the methods that  build upon them. They may  notably  lead to speed-ups when interfacing a neural network with a Kalman filter. 
\end{abstract}


\section{Introduction}
\label{sec:introduction}
The Kalman filter (KF) is a workhorse of state estimation for dynamical systems, which   is used in numerous technological fields.   Albeit optimal in the linear case, the KF critically relies on a large number of parameters, namely the covariance matrix of the initial error $P_0$, the covariance matrix $Q_n$  of the process noise and that of the measurement error $R_n$. Those parameters are generally unknown in practice, and need be estimated. Even when they are known through sensors' specifications, tuning them remains a challenge,   e.g., \cite{farrell2020imu}, and may require  manual ``tweaking" \cite{abbeel2005discriminative}.  Parameter inference is thus desirable and may then be achieved through  maximum likelihood (ML) or maximum a posteriori (MAP).   Historical approaches to maximum likelihood estimation of the Kalman  filter's parameters date back to the early days of Kalman filtering, and include descent optimization algorithms \cite{gupta1974computational} and derivative-free expectation-maximization approaches \cite{shumway1982approach}. 

Descent algorithms are based on the so-called sensitivity equations, first derived in  \cite{gupta1974computational}. They (forward) propagate recursively the derivative of all the filter’s outputs with respect to one scalar parameter. This allows for computing the sensitivity of   the likelihood  $\mathcal L$ of past estimates, with respect to parameter variations. More recent approaches focus on the numerical aspects and stability of KF implementation when computing the derivatives  \cite{kulikova2015constructing,tsyganova2017svd}. 

In deep learning,   backpropagation, sometimes called backprop, is the procedure for computing gradients  \cite{lecun2015deep}. Albeit often viewed as specific to neural networks, it is a way to apply the chain rule backwards, and can in fact  compute derivatives of any function. 
Using this approach, we  provide herein novel closed-form expression for gradient backpropogation through the Kalman filter. Those formulas show -  to the authors' surprise - that computing the sensitivity of the KF's estimates with respect to \emph{all} the  parameters at once, for instance matrix $\fracdlb{P_0}$,  is amenable to a  simple one-pass pipeline involving only matrix multiplications with matrices having the same dimension as the parameters. However, the proof is nontrivial. To actually derive closed-form expressions only, it   requires various manipulations of the KF equations based on  several  tricks, and judicious choices regarding the organization and the order in which  calculations are performed.  

Our contributions may be summarized as follows:
 \begin{itemize}
  \item We provide novel nontrivial equations in closed form for   computing the derivative backwards of a scalar  function $\Lc$ of a KF's trajectory of $N$ steps with respect to all its parameters. \item The derivative w.r.t. any  \emph{entire} parameter matrix   has computation cost $O(Nd^3)$ being similar to the KF's, where $d$ denotes the  dimension of the state space. 
  \item By contrast, previous classical sensitivity approaches, see e.g.,   \cite{gupta1974computational,kulikova2015constructing,tsyganova2017svd}, based on forward gradient computation have a cost $O(Nd^5)$ for the same derivative.
    \item  This is confirmed by a simple numerical experiment of dimension $d=6$:  the gradient of the  likelihood w.r.t. the  noise covariance matrix $R$ is speeded up by a   factor of 53 as compared to sensitivity equations.   
  \item In the same experiment we also achieve a speed-up   factor of 38 with respect to    state-of-the-art  numerical Pytorch's automatic differentiation (AD).
   \end{itemize}
   
   Beyond drastic computation speed-ups for sensitivity tuning, other motivations for using our technique are as follows: 
  \begin{itemize}
  \item Closed-form exact expressions allow practitioners for better control over the computational pipepline, while numerical methods (AD) are more ``black-box". 
    \item The numerical efficiency of our method may prove an enabler for fast backpropagation in  topical  applications that combine the KF with deep learning,   see \cite{haarnoja2016backprop,schymura2020dynbkf,chen2021dynanet,revach2022,
zhao2019deepmono,wu2021noisecov,brossard2020ai,jouaber2021nnakf}, 
 \end{itemize}

The paper is organized as follows. Preliminary notions are to be found in Section \ref{prelem:sec}.  
\Cref{main:sec} presents our main result and its proof is detailed in  the Appendix.   In \Cref{sec:experiments}, a numerical experiment on synthetic data illustrates the method.  

\section{Preliminaries}\label{prelem:sec}
 
 Before presenting our main result, we need a few primers. 

\subsection{Kalman filter equations}\label{sec:kf_eq}

For a dynamical system $  x_n= \Fn  x_{n-1} +\*B_n\*u_n+w_n$ with measurements $y_n =\Hn  x_n+v_n$, where $x_n$ is the state variable and $w_n,v_w$ white Gaussian noises, the state estimation $\xhpn$ of the KF is performed recursively through a prediction step
\begin{align}
\xhmn &= \Fn\xhpnn + \*B_n\*u_n\,,    \label{eq:innovx}\\    
\Pmn  &= \Fn\Ppnn \Fn^T + \Qn\,,     \label{eq:innovp}
\end{align}
followed by an update step
\begin{align}
\zn   & = \*y_n - \Hn\xhmn\,,   \label{eq:upz} \\
\Sn   & = \Hn\Pmn\Hn^T + \Rn\,,   \label{eq:ups} \\ 
\Kn   & = \Pmn\Hn^T\Sninv\,,   \label{eq:upk} \\
\xhpn & = \xhmn + \Kn\zn\,,   \label{eq:postx} \\
\Ppn  & = (\IKH)\Pmn\,,   \label{eq:upp}
\end{align}
where  $	\xhmn $  denotes the  predicted state based on past information $y_1,\dots,y_{n-1}$, with corresponding covariance matrix $
	\Pmn	$,  and  $
	\xhpn$ is the 	posterior state estimation of the KF in the light of latest measurement $y_n$, with corresponding covariance matrix $
	\Ppn	$. 
	$\Qn	$	is the covariance matrix of process noise $w_n$ and  $\Rn	$	that of  observation noise $v_n$ at time step $n$. The prediction error $\zn	$	is called  the innovation and has covariance matrix $\Sn$. $\Kn	$	is  the (optimal) Kalman gain.

\subsection{Loss function}
\label{sec:loss}
In this paper we consider a generic scalar cost  function defined over the entire past trajectory we call   loss, of the form 
\begin{align}
\Lc = \sum_{n=1}^{N} (\lnp  + \lnm),\label{loss:def}
\end{align}
where $\lnm$ and $\lnp$   respectively depend on the prior variables $\xhmn,\Pmn,\Rn,\yn$ and  the posterior variables $\xhpn,\Ppn$.  For parameter inference, the most widely considered loss  of the form \eqref{loss:def} is 
\begin{align}
\Linc = -\log{p (y_1, \ldots, y_N\mid\theta)} \, ,
\end{align}
which is the negative logarithm of the marginal likelihood (NLL) of the parameters encoded by the variable $\theta$, that may represent any tuning matrix    of the KF.   A standard calculation shows that up to a normalization constant it writes \begin{align}
\Linc = \sum_{n=1}^N \log\det S_n + z_n^T S_n^{-1} z_n \label{nrg:def} \, ,
\end{align}
 also known as the \emph{energy} function, see \cite{sarkka2013bayesian}, Chapter 12.  

Note  that, the interest of  the likelihood \eqref{nrg:def} is to allow for inference without a ground truth:   the performance criterion is a predictive performance on the observations $y_1,\dots,y_N$. However, in case of availability of a ground truth, that is, very accurate side information  provided by an extra sensor (only during the parameter inference  phase), our loss \eqref{loss:def} is versatile enough to encompass discrepancy of the state w.r.t. ground truth, as exploited in, e.g., \cite{abbeel2005discriminative,haarnoja2016backprop,brossard2020ai}.

\subsection{Matrix derivatives}
\label{sec:primer}
 
Let us consider a scalar function $ \mathcal L(M)$ depending on a $k\times m$  matrix  $M$ with $M=XY$ where $X$ and $Y$ are respectively $k\times d$ and $d\times m$ matrices.  We can define the gradient $\dv{\mathcal L}{M}$ as a $k\times m$ matrix, such that the entry $\bigl(\dv{\mathcal L}{M}\bigr)_{ij}$ is the partial derivative  of $ \mathcal L$ w.r.t. the element ${M}_{ij}$. Starting from  $\dv{\mathcal L}{M}$, we may compute $\dv{\mathcal L}{X}$   through the chain rule, which involves a large sum.  It turns out, though, that the matrix framework is convenient as   we may prove by inspection in the present case that $\dv{\mathcal L}{X} =\dv{\mathcal L}{M}Y^T$, that is, the derivative may be obtained through a simple matrix multiplication.

Indeed, the gradient $\dv{\mathcal L}{X}$ indicates how a small variation in the parameter $\delta X$ impacts the output. Using only matrix multiplications we have 
$\mathcal L(X+\delta X)\approx \mathcal L(X)+\mathrm{tr}\bigl((\fracd{\mathcal L}{X})^T\delta X\bigr)$, up to second order terms in  $\delta X$. This allows proving  general formulas for the derivative of the matrix composite function  $\mathcal L(M)$ with $M=\varphi(X)$,  see e.g., \cite{petersen2008matrix}.  
\begin{align}
	 	  M &=  X + Y \quad \Rightarrow \quad \hspace{.54cm} \fracdlb{X}=\fracdlb{M} \label{eq:xaddy} \\ 
	  	 M &=   XY \quad    \Rightarrow \quad  \hspace{.94cm}\fracdlb{X}= \fracdlb{M}Y^T \label{eq:xy} \\ 
	  M &=   YX \quad   \Rightarrow \quad  \hspace{.93cm} \fracdlb{X}=Y^T\fracdlb{M} \label{eq:yx} \\ 
	  M &=   YXY^T \quad   \Rightarrow \quad   \hspace{.45cm} \fracdlb{X}= Y^T\fracdlb{M}Y\label{eq:yxyt} \\ 
	 M &=   X^{-1} \quad  \Rightarrow  \quad  \hspace{.82cm} \fracdlb{X}= -M^{T} \fracdlb{M}M^{T}  \label{eq:invx}
 \\ 
	 \mathcal L&=\log\det M \quad  \hspace{.04cm}  \Rightarrow \quad \fracdlb{M}=M^{-T} \label{eq:logm}
	\end{align}
To be rigorous, one should write  for instance \eqref{eq:xy} as $\fracd{\!{L\circ\varphi}}{X}= \fracdlb{M}\Bigr|_{M=\varphi(X)}Y^T $. However, this would clutter notation.

 \subsection{Forward vs backward }

To differentiate the loss w.r.t. to a given state variable, covariance matrix, or noise parameter, we may distinguish between two approaches. A simple example is as follows:
\begin{align}
z=f(x_1,u),\quad  y=g(x_1,u) , \quad \Lc=h(y,z,u),
\label{smallex}\end{align}with $x_1,y,z,\Lc,$ scalars, and $u$ a parameter. Suppose we seek to compute    $\fracd{\Lc}{x_1}$, that is, to assess how a small variation  $\delta x_1$ generates a variation $\delta \Lc$.

\subsubsection{The forward method a.k.a sensitivity equations}

This method is the most intuitive and straightforward since it is a termwise differentiation of  the  equations. Consider   $x_1(\alpha)$  to allow for small changes in $x_1$. We may compute $\fracd{z}{\alpha}=\fracd{f}{x_1}\fracd{x_1}{\alpha}$ and $\fracd{y}{\alpha}=\fracd{g}{x_1}\fracd{x_1}{\alpha}$. A small variation may then be propagated forward, as  we have $\fracd{\Lc}{\alpha}=\fracd{h}{y}\fracd{y}{\alpha_i}+\fracd{h}{z}\fracd{z}{\alpha}$.  

When applied to the KF's equations \eqref{eq:innovx}-\eqref{eq:upp}   with particular  loss \eqref{nrg:def}, this leads to the  so-called sensitivity equations \cite{gupta1974computational}, see also \cite{sarkka2013bayesian}.  The problem with this method is that it does not scale well with the parameter's dimension, which may prove an impediment   when $x_1$ is a matrix    and $\Lc$ is a scalar.

 \subsubsection{The backward recursion method}
 Given  values of $x_1,u$, the values of $z$, $y$, and $\Lc$ are computed through \eqref{smallex}. Assume then that we want to assess how a small change in $x_1$ affects the computed value of $\Lc$. The chain rule may be computed backwards. Indeed, one may compute the gradients $\fracd{h}{y} $ and $\fracd{h}{z}$ at the obtained values for $y$ and $z$. We denote $\fracd{\Lc}{y}:= \fracd{h}{y}$ and $\fracd{\Lc}{z}:=\fracd{h}{z}$  and  store the corresponding values.    
  
 At this stage,  we have   independent ``contributions" of $\delta y$ and $\delta  z$ in the variation of $\Lc$. 
The connection appears when we propagate backwards. We may assess how a small variation $\delta x_1$ affects $\Lc$ through $z$ only, holding the other variables constant, as  expressed through the following notation
$$\bigl(\fracd{\Lc}{x_1}\bigr)^{(z)}=\fracd{\Lc}{z}\fracd{z}{x_1}\quad\text{(contribution via $z$)}\,,
$$using the already computed derivative $\fracd{\Lc}{z}$.  In the same way, 
$$\bigl(\fracd{\Lc}{x_1}\bigr)^{(y)}=\fracd{\Lc}{y}\fracd{y}{x_1}\quad\text{(contribution via $y$)}\,.$$The total variation of $\Lc$ in response to a variation of $x_1$ is obtained in turn by summing all those contributions, according to the chain rule, as \begin{align}\fracd{\Lc}{x_1}= \bigl(\fracd{\Lc}{x_1}\bigr)^{(y)}+\bigl(\fracd{\Lc}{x_1}\bigr)^{(z)}.\label{sum1:eq}
\end{align}

\subsection{Computation diagrams}

To  backward differentiate, dependency diagrams provide a useful guide, such as the one in \Cref{fig:graphformalism}  that  encapsulates  the dependencies in the KF equations \eqref{eq:innovx}-\eqref{eq:upp}, and on which we heavily rely  in the Appendix to prove the theorem. First, running the KF fixes all the values of the variables in the graph. To differentiate do as follows.  
Consider a function $\Lc$ of a given variable in the graph.  On may compute the value of $ \fracd{\Lc}{z_i}$, for all direct predecessors $z_i$ of this variable. For each predecessor $y_{ki}$ of $z_i$ one may in turn compute  $ \fracd{\Lc}{y_{ki}}=\fracd{\Lc}{z_{i}}\fracd{z_i}{y_{ki}}=\bigl(\fracd{\Lc}{y_{ki}}\bigr)^{(z_i)}$. This yields a front of values that back propagate.

To  differentiate with respect to any variable, say, $x$, we shall then identify all its direct successors $w_i$ in the graph, that is, all the arrows emanating from $x$. If all the   values  $\bigl(\fracd{\Lc}{x}\bigr)^{(w_i)}$ have been  obtained at previous backprop step,  we may write:
\begin{align}
\fracd{\Lc}{x}= \sum_{\text{$w_i$ direct successor of $x$}}\bigl(\fracd{\Lc}{x}\bigr)^{(w_i)}.\label{children:eq}
\end{align} 
This is exactly what has been done in \eqref{sum1:eq}.

\section{Main result}
\label{main:sec}

\begin{figure*}[!t]
\normalsize

\begin{align}
\intertext{\textbf{State variable  equations :}}
\fracdlb{\xhmn}  &= (\IKH)^T\fracdlb{\xhpn} + \fracd{\lnm}{\xhmn} \,,\label{Rn2001}\\
\fracdlb{\Pmn}   
				 &= (\IKH)^T \Big[\fracdlb{\Ppn} + \frac{1}{2} \fracdlb{\xhpn}\zn^T\Rn^{-1}\Hn + \frac{1}{2} \Hn^T\Rn^{-1}\zn(\fracdlb{\xhpn})^T  \Big](\IKH) + \fracd{\lnm}{\Pmn} \,, \label{Pmn_sym}\\
\fracdlb{\xhpnn} &=  \Fn^T \fracdlb{\xhmn}  + \fracd{l\ipostnn}{\xhpnn} \,, \label{eq1:eq} \\
\fracdlb{\Ppnn}  &=  \Fn^T \fracdlb{\Pmn} \Fn + \fracd{l\ipostnn}{\Ppnn} \,. \end{align} \begin{align}
\intertext{\textbf{Parameter   equations :}}
\fracdlb{\yn} &= \Kn^T\fracdlb{\xhpn} +\fracd{\lnm}{\yn}\, ,  \\
\fracdlb{\Qn} &= \fracdlb{\Pmn} \, ,  \\
\fracdlb{\Rn}   & =   
			 \Kn^T \fracdlb{\Ppn} \Kn  - \frac{1}{2} \Kn^T \fracdlb{\xhpn} \zn^T\Sn^{-1}  - \frac{1}{2}\Sn^{-1} \zn (\fracdlb{\xhpn})^T \Kn  + \fracd{\lnm}{\Rn}\,\label{Rn2000_sym22}  \\
\bigl(\fracdlb{\*L_n} &= 2\fracdlb{\Rn}\*L_n \quad\text{when optionally using the square-root form $\Rn = \*L_n\*L_n^T  $ } \bigr)\label{Rn2000_sym} 
\end{align} 

\hrulefill
\vspace*{4pt}
\end{figure*}
 
In this section we start with our main theorem, whose proof is postponed to the Appendix, and  discuss its application. 
\subsection{Main mathematical theorem}
\begin{theorem}[Backward parameter derivatives for the KF]\label{thm:bigthm}
Consider a loss function  $\Lc = \sum_{n=1}^{N} (\lnp  + \lnm) $,  where the functions involved are of the form  $\lnm(\xhmn,\Pmn,\Rn,\yn)$ and $\lnp(\xhpn,\Ppn)$. The matrix derivative of  $\Lc$  w.r.t. the matrix parameters at play  satisfy the following backward recursion. We start from $n=N$ and set  
$\fracdlb{\xhpn} =\fracd{l_{N|N}}{\xhpn} $ and $\fracdlb{\Ppn}=\fracd{l_{N|N}} {\Ppn}$,  and then go all the way to desired $n=k\leq N$ through  \eqref{Rn2001}-\eqref{Rn2000_sym}. In the case of static noise parameter, that is, $\forall n~R_n=R$, we may then apply the chain rule and readily compute $\fracdlb{R}=\sum_{n=1}^N \fracdlb{R_n}$, and  likewise $\fracdlb{Q}=\sum_{n=1}^N \fracdlb{Q_n}$.
\label{theorem}
\end{theorem} 
 
Note  \eqref{Rn2001}-\eqref{Rn2000_sym} give only the symmetric part of $\fracdlb{P},\fracdlb{R},\fracdlb{Q}$. Indeed, we want to retain the steepest descent direction, but within the space of symmetric matrices as the update must preserve the symmetry of the parameter matrices. We also note that to avoid storing $\Sn$ one can   replace $\Sn^{-1}$ with $\Rn^{-1}(I-\Hn\Kn)$ in \eqref{Rn2000_sym22}, as can be shown from Woodbury lemma.

\subsection{Numerical cost} 
The  equations provided by the theorem  are particularly economical computationally. 
Letting $d$ be the dimension of the state space, the numerical cost for computing all the derivatives is $O(Nd^3)$ with our method, since it involves only matrix multiplications at each of the $N$ steps. This is akin to  the cost of the KF Riccati equation, but even without matrix inversions.   This is in contrast with methods to date revolving around the well-known KF sensitivity equations \cite{gupta1974computational,sarkka2013bayesian}. Those methods compute partial derivatives w.r.t. the entries of the matrix parameters  {one by one}: to compute $\fracdlb{(P_0)_{ij}}$ one needs to implement a Riccati-like equation with computational cost   $O(Nd^3)$.  As there are $d(d+1)/2=O(d^2)$ entries the  numerical cost of computing $\fracdlb{P_0}$ for instance is $O(Nd^5)$.

 \subsection{Application to the energy  (NLL) function}

The theorem is versatile and accommodates a large range of loss functions, such as those based on  ground truth data, see \cite{schymura2020dynbkf,chen2021dynanet,revach2022,
zhao2019deepmono,wu2021noisecov,brossard2020ai,jouaber2021nnakf}. To apply it to   the energy function $\Linc$  defined by \eqref{nrg:def}, that reflects the negative log-likelihood  (NLL) of the parameters, we let $\lnp=0$, and  to $\lnm$ be  equal to $$ \log\det (\Hn\Pmn\Hn +\Rn)+ z_n^T (\Hn\Pmn\Hn +\Rn)^{-1} z_n,$$
with $\zn=\yn-\Hn\xhmn$. 
The derivative of the first term w.r.t. $\Rn$ is found to be  $\Sn^{-1}$ using \Cref{eq:xaddy}  and then the logdet derivative \eqref{eq:logm}.  As concerns the second term, let $w=z_n^T S_n^{-1} z_n$. We  have similarly $\frac{\partial w}{\partial\Rn} =\frac{\partial w}{\partial\Sn}=- \Sn^{-1}\fracd{w}{\Sn^{-1}}\Sn^{-1}=- \Sn^{-1}\zn\zn^T\Sn^{-1}$, proving
 \begin{align*}
 \fracd{\lnm}{\Rn}   = \Sn^{-1} - \Sn^{-1}\zn\zn^T\Sn^{-1} .
\end{align*}
Differentiating w.r.t. $(\Hn\Pmn\Hn)$ yields identical expressions, so $
 \fracd{\lnm}{\Pmn}  = \Hn^T \fracd{\lnm}{\Rn}  \Hn$. 
Completing those derivatives with standard gradient computation, we have
\begin{align}
& \fracd{\lnp}{\Ppn}  = 0, \qquad \fracd{\lnp}{\xhpn} = 0, \\
& \fracd{\lnm}{\Pmn}  = \Hn^T\Sn^{-1}\Hn - \Hn^T\Sn^{-1}\zn\zn^T\Sn^{-1}\Hn, \\
& \fracd{\lnm}{\xhmn} = -2\Hn^T\Sn^{-1}\zn, \\
& \fracd{\lnm}{\Rn}   = \Sn^{-1} - \Sn^{-1}\zn\zn^T\Sn^{-1}, \quad \fracd{\lnm}{\yn} = 2\Sn^{-1}\zn.
\end{align}
Substituting these into  the recursive equations of Theorem \ref{theorem}  we have all we need to computate   the derivative of the energy  (NLL) function  $\Linc$ of \Cref{nrg:def} with respect to all parameters of interest, including the measurements $\yn$.

Note that, we chose not to make the loss depend explicitly on $\Qn$, having no example in mind. However, $\!{L}$  depends on $\Qn$ indirectly, since $Q_n$ impacts   Kalman's estimates.

\subsection{Covariance matrices in square-root form}\label{num2:sec}
 
   When dealing  with positive semi-definite matrices, especially in the field of Kalman filtering, it is customary to work with square-root factors to enforce positive definiteness  \cite{bierman2006factorization}. Consider a positive semi-definite parameter matrix, say, $\*R$, and let us write it as $\*R = \*L\*L^T$ with $\*L$ a ``square-root factor''. 
   \begin{prop}The gradient w.r.t. $\*L$ is simply computed as 
\begin{align}
\fracd{\Lc}{\*L} = 2\fracd{\Lc}{\*R}\*L. \label{eq:grad_nnl}
\end{align}\end{prop}This justifies Eq. \eqref{Rn2000_sym}.

{\textbf{Proof}}: We recall that for a scalar function $f$ we have $f(R+\delta R)\approx f(R)+\mathrm{tr}\bigl((\fracd{f}{\*R})^T\delta R\bigr)$ neglecting second  order terms in $\delta R$ . Thus $f\bigl( (L+\delta L) (L+\delta L)^T\bigr)\approx f(LL^T)+\mathrm{tr}\bigl((\fracd{f}{\*R})^T(L\delta L^T+\delta L L^T)\bigr)= f(LL^T)+2\mathrm{tr}\bigl((\fracd{f}{\*R}L)^T\delta L\bigr)$, where we used  that we are interested only in the symmetric part of $\fracd{f}{\*R}$. Denoting $f(R)=\bar f(L)$ we  have $\bar f(L+\delta L)\approx\bar f(L)+  \mathrm{tr}\bigl((\fracd{\bar f}{\*L})^T \delta L\bigr) $, proving $\fracd{\bar f}{\*L}=2\fracd{f}{\*R}L$.

\section{Numerical experiments}
\label{sec:experiments}
  On a synthetic example, we compare our approach to:  \begin{itemize}\item the classical method of   sensitivity equations \cite{gupta1974computational,sarkka2013bayesian};  \item the state-of-the-art numerical method of   automatic differentiation (AD) using  the open framework PyTorch \cite{paszke2017automatic}. 
  \end{itemize}Although AD is being routinely used in machine learning to backpropagate in the KF, see e.g. \cite{schymura2020dynbkf,chen2021dynanet,revach2022,
zhao2019deepmono,wu2021noisecov,brossard2020ai,jouaber2021nnakf},  it turns out that using our analytical formulas significantly speeds up computations, and the gap with respect to the sensitivity equations method proves even bigger.

\subsection{Setting}

The experiment is performed on a simulated dataset  which provides a three-dimensional trajectory of a point in the 3D space with $N=1440$ time steps. 
The state consists of the six-dimensional vector  $x_n=(p_n,v_n)$ of position $p_n = (p_{n,x} \squad p_{n,y} \squad p_{n,z})$,  and velocity $v_n = (v_{n,x} \squad v_{n,y} \squad v_{n,z})$, and the known noisy accelerations are considered as inputs  $u_n = (a_{n,x} \squad a_{n,y} \squad a_{n,z})^T$.  
The dynamics write:
\begin{align}
p_n = p_{n-1} + \Delta t \, v_{n-1} + w_n^p \\
v_n = v_{n-1} + \Delta t \, u_{n-1} + w_n^v
\end{align} $\Delta t = 1 $ is the time step and $w_n$ is the process noise which is assumed Gaussian with zero mean and covariance  matrix $\*Q$.

The measurement $y_n$ is the noisy 3D position:
\begin{align}
\yn = p_n + \nu_n
\end{align}
where $\nu_n$ is the measurement noise which is assumed Gaussian with zero mean and unknown $3\times 3$ covariance matrix $\*R$. 

Our goal is to optimize the likelihood to estimate $\*R$.  We resort to the square-root form $\*R=\*L\*L^T$,   and compute the gradient w.r.t. ${\*L}$ of the energy (NLL) function \eqref{nrg:def}, that is,
\begin{align}
	\Linc = \sum_{n=1}^{N=1440} log|S_n| + z_n^T S_n^{-1}z_n  \, .
\end{align}
 To do so, we run a KF over the whole trajectory and we store the required values along the way depending on the method.
 AD is performed by PyTorch. For all methods we perform gradient descent $\*L \leftarrow\*L - \alpha\fracd{\Linc}{\*L}$, with a fixed step size $\alpha=0.005$. We use identical data, model  and filter, in order to obtain fair comparisons.

\subsection{Results}
\begin{figure}[h]
\centering
\includegraphics[scale=0.8]{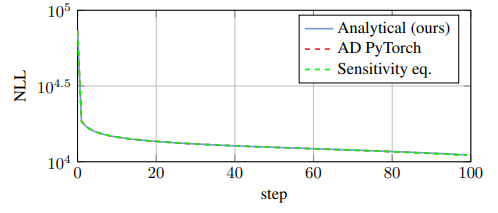}
\caption{Evolution at each gradient step of the NLL $\Lc$ in log scale using a) our analytical equations, b) PyTorch's Automatic Differentiation (AD), c) the sensitivity equations.   }
\label{fig:optimL}
\end{figure}
    
\begin{figure}[h]
\centering
\includegraphics[scale=.55]{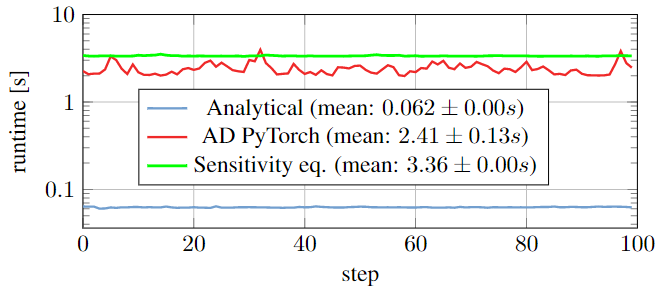}
\caption{Runtime in log scale at each gradient step with  our analytical method (blue line), AD PyTorch method (red line), and the sensitivity equations (green line). The analytical method is on average about 38 times faster than the AD method and 54 times faster than sensitivity equations. Moreover, AD shows a high variance of the  computation time whereas analytical methods require constant computation time.}
\label{fig:optimL_time}
\end{figure}

We compare the convergence of the loss $\Linc$ with respect to gradient descent iterations using our method, using sensitivity equations and using AD. The results are shown in \Cref{fig:optimL}. All yield identical gradients up to small numerical discrepancies, notably due to matrix inversion.    

The average computation time on CPU of each   step is displayed on   \Cref{fig:optimL_time}. The proposed method is respectively 38 times  and 53 times faster than AD and sensitivity method. This is consistent with our complexity analysis showing a gain by a $O(d^2)$ factor over sensitivity equation, with $d=6$. Further improvement    is thus anticipated for larger dimension $d$.  Moreover, the fluctuations of required process time for PyTorch backprop  indicate higher memory consumption. 

Finally, note we have implemented our algorithm in Python while coding in C++ may further enhance performance.
 
\section{Conclusion}

Using the powerful technique of backward matrix derivatives computation, and dependency graphs, we derived new closed-form expressions for parameter derivatives of a large class of loss functions for the Kalman filter, which bring benefits in terms of accuracy, stability of computation time, and above all  numerical cost.  

Fast computation time is critical in applications that combine KFs with neural networks. For instance, in   \cite{brossard2020ai}, the covariance matrix at each epoch consists of $\frac{15*16}{2}=120$ parameters and the trajectory is 1 hour long, which involves thousands of epochs. Hence, formulas yielding a significant reduction in training time may prove an enabling technology.

Analytic formulas also have the advantage over numerical differentiation to let the practioner keep control over what the method does. This may prove especially important  when the KF involves Cholesky or singular value decompositions (SVD) of the covariance matrix, see
\cite{bierman2006factorization}, which are ubiquitous in actual industrial implementations of the KF. With those decompositions the behaviour of AD seems harder to anticipate. 

 In the future, we would like to bring our exact backpropagation method to bear on the problem of training neural networks combined with Kalman filters in the difficult context of  high-dimensional sensors such as vision. The idea, pioneered in \cite{haarnoja2016backprop}, allows for the  best of both worlds: data-driven  feature learning combined with  optimal estimation and informed motion models.  In this regard, our contribution might play a key role for a new generation of learning-based  visual inertial systems (VINS),  see   \cite{zhao2019deepmono,bloesch2017iterated}.

  \appendix

  \section{Appendix: Derivation of the  main result}
  \label{sec:proof}
Let us  go through the  steps needed to prove Theorem \ref{thm:bigthm}. 
\subsection{Initialization} 

 The recursion is initialized by computing the derivatives \eqref{init:ref} at time $n=N$, that is, by differentiating $l_{N|N}$ (only). 
 Then, at each step $1\leq n\leq N$  of the backward recursion,  we assume known the following
\begin{align}
\fracdlb{\xhpn}  ,\qquad \fracdlb{\Ppn} \,, \label{init:ref}
\end{align} 
which  reflect how a small variation in $\xhpn$ and $\Ppn$ affects the loss $\Lc$, ignoring variables up to previous step $n-1$. We then  evaluate  all the desired derivatives at previous step. 
 
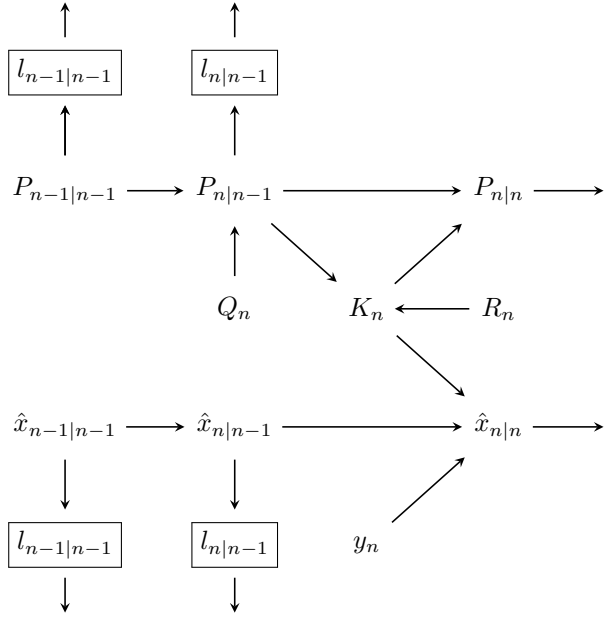
\begin{figure}
\begin{center}
    \begin{tikzpicture}
	[
	nostyle/.style={},
	squarednode/.style={rectangle, draw=red!60, fill=red!5, very thick, minimum size=5mm},
	            > = stealth,
	            shorten > = 1pt, 
	            auto,
	            node distance = 0.7cm, 
	            minimum height = 0.85cm,
	            semithick 
	]

	\node[minimum height = 0cm]		(11)					  {};
	\node[nostyle]		(21)		[below=0.5cm of 11] {$\boxed{l\ipostnn}$};
	\node[nostyle]		(31)		[below=of 21] {$\*P\ipostnn$};
	\node[nostyle]		(41)		[below=of 31] {};
	\node[nostyle]		(51)		[below=of 41] {$\xhpnn$};
	\node[nostyle]		(61)		[below=of 51] {$\boxed{l\ipostnn}$};
	\node[minimum height = 0cm]		(71)		[below=0.5cm of 61] {};
	
	\node[minimum height = 0cm]		(12)		[right=2cm of 11] {};
	\node[nostyle]		(22)		[below=0.5cm of 12] {$\boxed{\lnm}$};
	\node[nostyle]		(32)		[below=of 22] {$\Pmn$};
	\node[nostyle]		(42)		[below=of 32] {$\Qn$};
	\node[nostyle]		(52)		[below=of 42] {$\xhmn$};
	\node[nostyle]		(62)		[below=of 52] {$\boxed{\lnm}$};
	\node[minimum height = 0cm]		(72)		[below=0.5cm of 62] {};
	
	\node[minimum height = 0cm]		(13)		[right=1.5cm of 12] {};
	\node[nostyle]		(23)		[below=0.5cm of 13] {};
	\node[nostyle]		(33)		[below=of 23] {};
	\node[nostyle]		(43)		[below=of 33] {$\Kn$};
	\node[nostyle]		(53)		[below=of 43] {};
	\node[nostyle]		(63)		[below=of 53] {$\yn$};
	\node[minimum height = 0cm]		(73)		[below=0.5cm of 63] {};
	
	\node[minimum height = 0cm]		(14)		[right=1.5cm of 13] {};
	\node[nostyle]		(24)		[below=0.5cm of 14] {};
	\node[nostyle]		(34)		[below=of 24] {$\Ppn$};
	\node[nostyle]		(44)		[below=of 34] {$\Rn$ };
	\node[nostyle]		(54)		[below=of 44] {$\xhpn$};
	\node[nostyle]		(64)		[below=of 54] {};
	\node[minimum height = 0cm]		(74)		[below=0.5cm of 64] {};
	
	\node[minimum height = 0cm]		(15)		[right=1.3cm of 14] {};
	\node[nostyle]		(25)		[below=0.5cm of 15] {};
	\node[nostyle]		(35)		[below=of 25] {};
	\node[nostyle]		(45)		[below=of 35] {};
	\node[nostyle]		(55)		[below=of 45] {};
	\node[nostyle]		(65)		[below=of 55] {};
	\node[minimum height = 0cm]		(75)		[below=of 65] {};
	
	\path[<-] (11) edge node {} (21);
	\path[<-] (21) edge node {} (31);
	\path[->] (51) edge node {} (61);
	\path[->] (61) edge node {} (71);
	
	\path[<-] (12) edge node {} (22);
	\path[<-] (22) edge node {} (32);
	\path[->] (42) edge node {} (32);
	\path[->] (52) edge node {} (62);
	\path[->] (62) edge node {} (72);
	
	\path[->] (31) edge node {} (32);
	\path[->] (51) edge node {} (52);
	
	\path[->] (32) edge node {} (43);
	\path[<-] (21) edge node {} (31);
	
	\path[->] (32) edge node {} (34);
	\path[<-] (43) edge node {} (44);
	\path[<-] (54) edge node {} (43);
	\path[<-] (54) edge node {} (63);
	\path[->] (52) edge node {} (54);

	\path[->] (34) edge node {} (35);
	\path[->] (54) edge node {} (55);
	\path[->] (43) edge node {} (34);
	
	\end{tikzpicture}
	\caption{Computation diagram for the KF's recursive equations \eqref{eq:innovx}-\eqref{eq:upp}. All arrows pointing to a variable indicate the variables explicitly involved in its expression. When a given variable (a node) slightly varies - let us call it $x$ -  it affects all the variables to which an arrow emanating from $x$ points (the direct successors). The variation of $x$ then contributes to a variation of  $\Lc$ through all its direct successors, by adding all the corresponding contributions, as indicated by the chain rule, see \eqref{children:eq}.} \label{fig:graphformalism}
 \end{center}
\end{figure}

\subsection{Derivatives  w.r.t.   $\yn$  and $\xhmn$} Looking at the diagram of Figure 
\ref{fig:graphformalism} we see  $\yn$ has  $\xhpn$  as sole successor, and  $\xhmn$ additionally has $\lnm$.   
$\xhpn$ is a function of those variables via   update equation \eqref{eq:postx}:
\begin{align}
\xhpn   =  \xhmn + \Kn\zn=(\Id-\Kn \Hn ) \xhmn+ \Kn \yn. \label{eq:innovdev}
\end{align}
$\yn$ acts on $\xhpn$ via the product $\Kn\yn$. From  \eqref{eq:yx} we get
\begin{align}
\fracdlb{\yn}= \Kn^T\fracdlb{\xhpn}. 
\end{align}
Using   \eqref{eq:innovdev} and 
\eqref{eq:yx}  yields $\bigr(\fracdlb{\xhmn}\bigr)^{(\xhpn)} = (\Id-\Kn\Hn)^T\fracdlb{\xhpn} $, where we recall the notation indicates a contribution via $\xhpn $.  Besides,     there is one  other arrow emanating from $\xhmn$, which points to $\lnm$. Summing those contributions according to the rule \eqref{children:eq}  yields:
\begin{align}
\fracdlb{\xhmn} = (\Id-\Kn\Hn)^T\fracdlb{\xhpn} + \fracd{\lnm}{\xhmn}\,. 
\end{align}

\subsection{Derivative w.r.t. $\Pmn$} 

Looking at the  successors of $\Pmn$ on the  diagram of  \Cref{fig:graphformalism}, we see a small variation in  $\Pmn$  generates a variation in $\Lc$ through $\Ppn$,   through  $\xhpn$, whose  calculation depends on the gain $\Kn$ as a function of $\Pmn$, and directly through $\lnm$. 
The relevant update equations of the KF  which encapsulate these dependancies are \eqref{eq:ups}-\eqref{eq:upp}.

\subsubsection{Contribution through $\Ppn$}
 
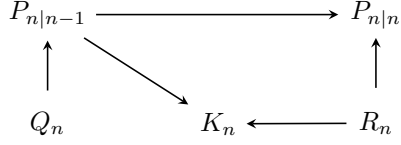
\begin{figure}
\begin{center}
    
	\begin{tikzpicture}
	[
	nostyle/.style={},
	squarednode/.style={rectangle, draw=red!60, fill=red!5, very thick, minimum size=5mm},
	            > = stealth, 
	            shorten > = 1pt, 
	            auto,
	            node distance = 0.7cm, 
	            minimum height = 0.85cm,
	            semithick 
	]	 
	
	\node[minimum height = 0cm]		(32)		  {$\Pmn$};
	\node[nostyle]		(42)		[below=of 32] {$\Qn$};
	 
	 	\path[->] (42) edge node {} (32);

		\node[minimum height = 0cm]		(13)		[right=1.5cm of 32] {};
	\node[nostyle]		(43)		[below=0.89cm of 13]  {$\Kn$};

	\path[->] (32) edge node {} (43);
	
	\node[minimum height = 0cm]		(34)		[right=1.5cm of 13] {$\Ppn$};
	\node[nostyle]		(44)		[below=of 34] {$\Rn$ };

	\path[->] (32) edge node {} (34);
	\path[<-] (43) edge node {} (44);
	\path[<-] (34) edge node {} (44);

	\end{tikzpicture}
	\caption{Zoom on the part of the diagram modified   after \Cref{eq:upp} has been transformed into \Cref{eq:inf_pn}. This allows for $\Ppn$ to be unaffected by a variation of variable $\Kn$, rendering them ``independent''.} \label{fig:graphformalism2}
 \end{center}

\end{figure}

Regarding the direct dependence of $\Ppn$ on $\Pmn$,   using \Cref{eq:upp} may prove convoluted as $\Kn$ also nontrivially depends on $\Pmn$. A   trick to get around this issue consists in resorting to the information form of the KF equations instead:
\begin{align}
(\Ppn)^{-1} &= (\Pmn)^{-1} + \Hn^T\Rn^{-1}\Hn\,. \label{eq:inf_pn}
\end{align}
This allows modifying the dependencies in the original diagram of \Cref{fig:graphformalism} as illustrated by the diagram of  \Cref{fig:graphformalism2}. 

Starting from \Cref{eq:inf_pn}, and using  the matrix derivative formula \eqref{eq:xaddy}, proves
\begin{align}
\bigl(\fracdlb{\Pmn^{-1}}\bigr)^{(\Ppn)} & =  \fracdlb{\Ppn^{-1}},  \label{eq:gradPm}
\end{align}  
where the notation indicates we evaluate the contribution of a variation of $\Pmn$ on the loss  through the node $\Ppn$. Using the matrix   formula \eqref{eq:invx} we then get
 \begin{equation}\begin{aligned}\bigl(\fracdlb{\Pmn}\bigr)^{(\Ppn)}&=- \Pmn^{-1} \bigl(\fracdlb{\Pmn^{-1}}\bigr)^{(\Ppn)} \Pmn^{-1}  \,,\\&=\Pmn^{-1} \Ppn  \fracdlb{\Ppn}  \Ppn \Pmn^{-1}\\&=(\IKH)^T  \fracdlb{\Ppn}   (\IKH). \end{aligned} \label{eq:dp}\end{equation}

\subsubsection{Contribution through $\xhpn$}
Let us deal with the other arrow emanating from $\Pmn$. A variation in  $\Pmn$  also  affects  $\xhpn$ through $\Kn$, see \eqref{eq:innovdev}.  The problem is that the dependence of $\Kn$ on $\Pmn$ is complicated, see  \eqref{eq:ups}-\eqref{eq:upk}. Another trick may be used in the way we just did. Indeed,  the gain alternatively writes $\Kn=\Ppn \Hn^T \Rn^{-1}$, see e.g., \cite{deyst1968conditions}.  This allows for rewriting \eqref{eq:innovdev} as	\begin{align}\xhpn  = \xhmn + \Ppn \Hn^T \Rn^{-1}\zn.\label{eq49:eq}\end{align}Letting $u=\Ppn \Hn^T \Rn^{-1}\zn$, we get    $\bigl(\fracdlb{\Ppn}\bigr)^{(\xhpn)}  \stackrel{\eqref{eq:xy}}{=}  \bigl(\fracdlb{u}\bigr)^{(\xhpn)} \zn^T\Rn^{-1}\Hn
\stackrel{\eqref{eq49:eq},\eqref{eq:xaddy}}{=}  \bigl(\fracdlb{\xhpn}\bigr)^{(\xhpn)} \zn^T\Rn^{-1}\Hn$. Finally, using that $\Pmn$ and $\Ppn$ are related by \eqref{eq:inf_pn} and redoing calculation \eqref{eq:dp} yields
	$ \bigl(\fracdlb{\Pmn}\bigr)^{(\xhpn)} =(\IKH)^T\fracdlb{\xhpn} \zn^T\Rn^{-1}\Hn(\IKH).$

The rule \eqref{children:eq}  applied to all   the successors of $\Pmn$ yields
\begin{align}
 &\fracdlb{\Pmn} = (\IKH)^T \bigg[ \fracdlb{\Ppn} + \fracdlb{\xhpn} \zn^T\Rn^{-1}\Hn \bigg] \nonumber\\ &\quad 
	(\IKH) + \fracd{\lnm}{\Pmn} \, , 
	\label{eq:Pmn_asym}
     	\end{align}
	where we factorized by $(\IKH)^T$ on the left and $(\IKH)$ on the right. 
We note  the expression $\fracdlb{\xhpn} \zn^T\Rn^{-1}\Hn$ in \eqref{eq:Pmn_asym} is not symmetric. We have performed an unconstrained calculation but parameter $\Pmn$ is symmetric. Thus, we retain only    the symmetrical part, so that a gradient  shall preserve symmetry of  $\Pmn$, leading to \eqref{Pmn_sym}. 

 \subsection{Derivative w.r.t. parameter matrix $\Rn$}
On  Fig. \ref{fig:graphformalism2} we see the successors of $\Rn$ are  $\Ppn$,  and $\Kn$ (which points to $\xhpn$ in turn, see Figure  \ref{fig:graphformalism}). 
Regarding $\Ppn$, using  formulas \eqref{eq:xaddy}, \eqref{eq:yxyt},   \eqref{eq:invx},   \Cref{eq:inf_pn} yields $\fracdlb{ \Rn^{-1}}=\Hn \fracdlb{\Ppn^{-1} } \Hn^T=-\Hn \Ppn \fracdlb{\Ppn } \Ppn\Hn^T$ so that from  \eqref{eq:invx} we have $\bigl(\fracdlb{\Rn}\bigr)^{(\Ppn)}=\Rn^{-1} \Hn \Ppn \fracdlb{\Ppn} \Ppn \Hn^T  \Rn^{-1}=\Kn^T \fracdlb{\Ppn} \Kn.$

For successor $\xhpn$ via $\Kn$, as $\Ppn$ depends also on $\Rn$ in \eqref{eq49:eq}, we'd rather get back to \eqref{eq:innovdev} to derivate w.r.t. $\Rn$. 
 We have
$\bigl(\fracdlb{\Rn}\bigr)^{(\xhpn)}\stackrel{\eqref{eq:ups}}{=}\fracdlb{\Sn}\stackrel{\eqref{eq:invx}}{=}-\Sn^{-1} \fracdlb{\Sn^{-1}}\Sn^{-1}\stackrel{\eqref{eq:upk}}{=}-\Sn^{-1}\Hn \Pmn\fracdlb{\Kn}\Sn^{-1}
=-\Kn^T \fracdlb{\Kn}\Sn^{-1}
\stackrel{\eqref{eq:gradK}}{=}-\Kn^T\fracdlb{\xhpn}\zn^T\Sn^{-1}.$ The latter stems from the fact a variation in $\Kn$ impacts $\mathcal L$ via node $\xhpn$ using \eqref{eq:innovdev}. It yields
\begin{align}
\fracdlb{\Kn} = \fracdlb{\xhpn}\zn^T\,. \label{eq:gradK}
\end{align}
Summing the contributions yields: $
\bigl(\fracdlb{\Rn}\bigr)^{(\Ppn)}+\bigl(\fracdlb{\Rn}\bigr)^{(\xhpn)}  
	=	\Kn^T \fracdlb{\Ppn} \Kn - \Kn^T \fracdlb{\xhpn} \zn^T\Sn^{-1} $. By additionally accounting for the explicit dependance of $\lnm$ on $\Rn$, and applying our rule \eqref{children:eq}, we get the desired equation. Finally, since $\Rn$ is symmetric, we symmetrize as before, leading to   \Cref{Rn2000_sym22}.

\subsection{Derivative w.r.t. $\xhpnn$, $\Ppnn$ and $\Qn$}

The diagram of \Cref{fig:graphformalism}   shows the simple dependencies encapsulated in the KF propagation equations \eqref{eq:innovx}-\eqref{eq:innovp}. Not omitting the direct dependency on $\xhpnn $, \eqref{eq:innovx} yields
\begin{align}
\fracdlb{\xhpnn} &= \Fn^T \fracdlb{\xhmn} + \fracd{l_{n-1|n-1}}{\xhpnn}\, ,
\end{align}and the counterpart  for \eqref{eq:innovp} writes using \eqref{eq:yxyt}:\begin{align}
\fracdlb{\Ppnn} &= \Fn^T \fracdlb{\Pmn} \Fn + \fracd{l_{n-1|n-1}}{\Ppnn} \, .
\end{align}
Finally using matrix formula \eqref{eq:xaddy} and \Cref{eq:innovp}  merely  yields $
	\fracdlb{\Qn} = \fracdlb{\Pmn}$. This completes the proof.

\bibliographystyle{unsrtnat}

\end{document}